\documentclass[10 pt]{article}

\usepackage{amsmath, amssymb, amsfonts, amsthm}
\usepackage{graphics}
\usepackage[margin=1in]{geometry}
\usepackage{epstopdf}
\usepackage{setspace}
\usepackage{fancyhdr}
\usepackage{rotating}
\usepackage{multirow}
\usepackage{graphicx}
\usepackage{algorithm}
\usepackage[margin=0.5cm]{subcaption}
\usepackage[margin=1cm]{caption}
\newcommand{\argmin}{\operatornamewithlimits{argmin}}

\title{Separation of Undersampled Composite Signals using the Dantzig Selector with Overcomplete Dictionaries\thanks{This paper is a preprint of a paper accepted by IET Signal Processing and is subject to Institution of Engineering and Technology Copyright. When the final version is published, the copy of record will be available at IET Digital Library.  Cleared for public release by WPAFB Public Affairs on 28 Aug 14.  Case Number: 88ABW-2014-4075.
This research is supported in part by an award from the National Research Council via the Air Force Office of Scientific Research and by the US National Science Foundation under grant DMS-1115523.}}

\author{
Ashley Prater\thanks{Correspondence Author.  Air Force Research Laboratory, High Performance Systems Branch, 525 Brooks Rd, Rome, NY 13441. \tt{ashley.prater.3@us.af.mil}} \and Lixin Shen\thanks{Department of Mathematics, Syracuse
University, 215 Carnegie Building, Syracuse, NY 13244. \tt{lshen03@syr.edu} }
}

\begin{document}
\maketitle

\begin{abstract}
In many applications one may acquire a composition of several signals that may be corrupted by noise, and it is a challenging problem to reliably separate the components from one another without sacrificing significant details.   Adding to the challenge, in a compressive sensing framework, one is given only an undersampled set of linear projections of the composite signal.  In this paper, we propose using the Dantzig selector model incorporating an overcomplete dictionary to separate a noisy undersampled collection of composite signals, and present an algorithm to efficiently solve the model.

The Dantzig selector is a statistical approach to finding a solution to a noisy linear regression problem by minimizing the $\ell_1$ norm of candidate coefficient vectors while constraining the scope of the residuals.  If the underlying coefficient vector is sparse, then the Dantzig selector performs well in the recovery and separation of the unknown composite signal.   In the following, we propose a proximity operator based algorithm to recover and separate unknown noisy undersampled composite signals through the Dantzig selector.  We present numerical simulations comparing the proposed algorithm with the competing Alternating Direction Method, and the proposed algorithm is found to be faster, while producing similar quality results.  Additionally, we demonstrate the utility of the proposed algorithm in several experiments by applying it in various domain applications including the recovery of complex-valued coefficient vectors, the removal of impulse noise from smooth signals, and the separation and classification of a composition of handwritten digits.

\end{abstract}

\section{Introduction}
This paper considers the problem of separating a composite signal through the recovery of an underlying sparse coefficient vector by using the Dantzig selector given only an incomplete set of noisy linear random projections.  That is, we discuss the estimation of a coefficient vector $c \in \mathbb{C}^{q}$ given the vector
\begin{equation}\label{eqn:linear}
	y = X\beta + z,
\end{equation}
where $X\in\mathbb{R}^{n\times p}$ is a sensing matrix with $n\leq p$, $z$ is a collection of i.i.d.\ $\sim N(0,\sigma^2)$ random variables, and the unknown signal $\beta \in \mathbb{R}^p$ admits the sparse representation $\beta = Bc$ for a known overcomplete dictionary $B \in \mathbb{C}^{p \times q}$. The individual signals composed to form $\beta$ can then be recovered from $c$ and $B$.  Since Equation~\eqref{eqn:linear} is underdetermined yet consistent, it presents infinitely many candidate signals $\beta$ and coefficient vectors $c$. 

The Dantzig selector was introduced in \cite{C05thedantzig} as a method for estimating a sparse parameter $\beta \in \mathbb{R}^{p}$ satisfying~\eqref{eqn:linear}. Discussions on the Dantzig selector, including comparisons to the least absolute shrinkage and selection operator (LASSO), can be found in~\cite{bickelDiscussion, caiDiscussion, C05thedantzig, candesRejoinder, efronDiscussion, friedlanderDiscussion, meinsDiscussion, ritovDiscussion}. 
Both the Dantzig selector and LASSO aim for sparse solutions, but whereas LASSO tries to match the image of candidate vectors close to the observations, the Dantzig selector aims to bound the predictor of the residuals. When tuning parameters in LASSO and the Dantzig selector model are set properly, the LASSO estimate is always a feasible solution to the Dantzig selector minimization problem, although it may not be an optimal solution. Furthermore, when the corresponding solutions are not identical, the Dantzig selector solution is sparser than the LASSO solution in terms of the $\ell_1$ norm \cite{James-Radchenko-Lv:jrss:09}. Recently, the Dantzig selector model has been applied for gene selection in cancer classification \cite{Zheng-Liu:CBM:11}. 

Classical compressive sensing theory guarantees the recovery of a sparse signal given only a very small number of linear projections under certain conditions \cite{Baraniuk07compressivesensing, linearprogramming, candesOptimal, Donoho06stablerecovery}. However, very seldomly is a naturally encountered signal perfectly sparse relative to a single basis.  Therefore, a number of works have considered the recoverability of signals that are sparse relative to an overcomplete dictionary that is formed by the concatenation of several bases or Parseval frames \cite{candesRedundant, Donoho_eladm, Donoho06stablerecovery, Elad2002generalized, kutyniok2011, rauhut}.   In this work, we propose and analyze a Dantzig selector model inspired by the above applications of overcomplete dictionaries in compressive sensing, and develop an algorithm for finding solutions to this model.

The following notation will be used.  The absolute value of a scalar $\alpha$ is denoted by $|\alpha|$, and the number of elements in a set $T$ is denoted by $|T|$.  The smallest integer larger than the real number $\alpha$ is denoted by $\lceil \alpha \rceil$.  The $i^{\text{th}}$ element of a vector $x$ is denoted by $x(i)$, and the $i^{\text{th}}$ column of a matrix $A$ is denoted by $A_i$.  The support of a vector $x$ is given by $\text{supp}(x) = \{i: x(i) \neq 0\}$.  The $\ell_1$, and $\ell_2$ vector norms, denoted by $\|\cdot\|_1$, and $\|\cdot\|_2$ respectively,  are defined by
\[ \|x\|_1 = \sum_{i=1}^n \left|x(i)\right|, \; \|x\|_2 = \left(\sum_{i=1}^n \left| x(i)\right|^2\right)^{\frac{1}{2}}, \]
for any vector $x \in \mathbb{C}^n$.  For matrices $A, B$ with the same number of rows, $\begin{bmatrix} A & B\end{bmatrix}$ is the horizontal concatenation of $A$ and $B$.  Similarly, $\begin{bmatrix} A \\ B\end{bmatrix}$  is the vertical concatenation of $A$ and $B$, provided each has the same number of columns. The conjugate transpose of a matrix $A$ is denoted by $A^{\top}$.

The rest of the paper is organized as follows.   In Section 2, the Dantzig selector model incorporating overcomplete dictionaries is introduced.  In Section 3, we present an algorithm used to find solutions to the proposed model.  Section 4 presents several numerical experiments demonstrating the appropriateness of the model and the accuracy of the results produced by the presented algorithm.  In simulations using real-valued matrices in the overcomplete dictionary, we compare the efficiency and accuracy of the presented method with the competing Alternating Direction Method.  Additionally, we demonstrate the utility of the proposed algorithm in several experiments by applying it in various domain applications including the recovery of complex-valued coefficient vectors, the removal of impulse noise from smooth signals, and the separation and classification of a composition of handwritten digits.  We close the paper with some remarks and possible future directions.

\section{The Dantzig selector model incorporating overcomplete dictionaries}
In this section, we present a Dantzig selector model incorporating overcomplete dictionaries that can be used to recover an unknown signal and reliably separate overlapping signals.

Suppose the unknown composite signal $\beta$ is measured via $y = X\beta+z$, where $X$ is an $n\times p$ sensing matrix and $z$ models sensor noise, and suppose an overcomplete dictionary $B$ is known such that $\beta = Bc$ for some sparse $c$.  Although $\beta$ and $c$ are not known, it is reasonable in many applications to know or suspect the correct dictionary components.  For example, if the signals of interest appear to be sinusoids with occasional spikes as in Figure~\ref{fig:experiment 3}, one should use a dictionary that is a concatenation of a discrete Fourier transform component and a standard Euclidean basis component.  In the following, let $q = 2p$ and assume the $p\times q$ dictionary $B$ is formed by a horizontal concatenation of a pair of orthonormal bases, $B = \begin{bmatrix}\Phi & \Psi\end{bmatrix}$, and the components of $\beta$ admit the sparse representations $\beta_{\Phi} = \Phi c_{\Phi}$ and $\beta_{\Psi} = \Psi c_{\Psi}$, with $\beta = \beta_{\Phi} + \beta_{\Psi}$ and $c = \begin{bmatrix} c_{\Phi}^\top & c_{\Psi}^\top \end{bmatrix}^\top$.  More succinctly,
$$
\beta = \begin{bmatrix} \Phi & \Psi \end{bmatrix} \begin{bmatrix} c_{\Phi} \\ c_{\Psi} \end{bmatrix}.
$$

To recover $c$, and therefore also $\beta$ and the components $\beta_{\Phi}$ and $\beta_{\Psi}$, from the observations $y$, we propose using a solution to the Dantzig selector model (see~\cite{C05thedantzig}) with an overcomplete dictionary
\begin{equation}\label{model:dantzig overcomplete}
	c \in \min_{c\in\mathbb{C}^{2p}} \left\{ \|c\|_1 : \| D^{-1} B^\top X^\top\left(XBc-y\right)\|_{\infty}\leq \delta \right\},
\end{equation}
where the diagonal matrix $D\in\mathbb{R}^{q\times q}$ with entries $d_{jj} = \text{diag}\{\|(XB)_j\|_2\}$ normalizes the sensing-dictionary pair.  Although Model~\eqref{model:dantzig overcomplete} is expressed using an overcomplete dictionary with two representation systems, one could generalize the model to accomodate more systems.

If the elements of $X$ are independent and identically distributed random variables from a Gaussian or Bernoulli distribution, and $B$ contains elements of fixed, nonrandom bases, then $D$ is invertible.  To see this, note that $d_{jj} = 0$ if and only if $\langle (X^\top)_i, B_j \rangle = 0$ for all  $i \in \{1,2,\ldots,n\}$.   However, since a random sensing matrix is largely incoherent with, yet not orthogonal to any fixed basis~\cite{candesRandom, donohohuo, greed}, it follows that $d_{jj}\neq 0$ for each $j$, ensuring $D$ is invertible.  Employing a sensing matrix whose entries are i.i.d.\ random variables sampled from a Gaussian or Bernoulli distribution, paired with an overcomplete dictionary formed by several bases or parseval frames has the added benefit of giving small restricted isometry constants, which in turn improves the probability of successful recovery of the coefficient vector via $\ell_1$ minimization.  More on these concepts, now standard in compressive sensing literature, can be found in~\cite{abol, Baraniuk07compressivesensing, candesOptimal,  candesRedundant, Donoho_eladm, Donoho06stablerecovery, donohohuo, Elad2002generalized, kutyniok2011, rauhut}.


\section{A proximity operator based algorithm}
To compute the Dantzig selector, we characterize a solution of Model~\eqref{model:dantzig overcomplete} using the fixed point of a system of equations involving applications of the proximity operator to the $\ell_1$ norm.   In this section we describe the system of equations and their relationship to the solution of Model~\eqref{model:dantzig overcomplete} and present an algorithm with an iterative approach for finding these solutions.

Let $A = D^{-1} B^\top X^\top XB$, and define the vector $\gamma = D^{-1} B^\top X^\top y$ and the set $\mathcal{F} = \{ c: \|c-\gamma\|_{\infty}\leq \delta\}$. The indicator function $\iota_{\mathcal{F}}:\mathbb{C}^{2p}\rightarrow \{0,\infty\}$ is defined by
\[
	\iota_{\mathcal{F}}(u) = \begin{cases}0, &\text{if } u \in \mathcal{F} \\ +\infty, &\text{if } u\notin \mathcal{F} \end{cases}
\]
and the proximity operator of a lower semicontinuous convex function $f$ with parameter $\lambda\neq 0$ is defined by
\[ \text{prox}_{\lambda f}(x) = \argmin_{u \in \mathbb{C}^{2p}} \left\{ \frac{1}{2\lambda}\|u - x\|^2_2 + f(u) \right\}. \]
Then Model~\eqref{model:dantzig overcomplete} can be expressed in terms of the indicator function as
\begin{equation}\label{eqn:indicator function model}
	c \in \min_{c\in\mathbb{C}^{2p}}\left\{ \|c\|_1 + \iota_{\mathcal{F}}(Ac) \right\}.
\end{equation}
If $c$ is a solution to Model~\eqref{eqn:indicator function model}, then for any $\alpha, \lambda>0$ there exists a vector $\tau \in \mathbb{C}^{2p}$ such that
\[
	c = \text{prox}_{\frac{1}{\alpha}\|\cdot\|_1}\left( c - \frac{\lambda}{\alpha}A^\top\tau\right) \text{ and } \tau = \left(I - \text{prox}_{\iota_{\mathcal{F}}}\right)\left(Ac + \tau\right).
\]
Furthermore, given $\alpha$ and $\lambda$, if $c$ and $\tau$ satisfying the above equations exist, then $c$ is a solution to~\eqref{eqn:indicator function model}, and therefore also to~\eqref{model:dantzig overcomplete}.  Using the fixed-point characterization above, the $(k+1)^{\text{th}}$ iteration of the proximity operator algorithm to find the solution of the Dantzig selector model incorporating an overcomplete dictionary is
\begin{equation}\label{eq:itr}
	\begin{cases}
	&c^{k+1} = \text{prox}_{\frac{1}{\alpha}\|\cdot\|_1}\left( c^k - \frac{\lambda}{\alpha}A^\top(2\tau^k - \tau^{k-1})\right), \\ & \tau^{k+1} = \left(I - \text{prox}_{\iota_{\mathcal{F}}}\right)\left(Ac^{k+1} + \tau^k\right). \end{cases}
\end{equation}
If $\lambda/\alpha < 1/\|A\|_2^2$, the sequence $\{(c^k,\tau^k)\}$ converges.  The proof follows those in~\cite{chambolle, micchelli}. We remark that the proximity operators appearing in Equation~\eqref{eq:itr} can be efficiently computed. More precisely, for any positive number $\lambda$ and any vector $u \in \mathbb{C}^d$,
$$
\text{prox}_{\lambda\|\cdot\|_1} (u) = \begin{bmatrix}\text{prox}_{\lambda |\cdot|} (u_1) &  \text{prox}_{\lambda |\cdot|} (u_2) & \cdots &\text{prox}_{\lambda |\cdot|} (u_{2p}) \end{bmatrix}^\top,
$$
and
$$
\text{prox}_{\iota_{\mathcal{F}}} (u) = \begin{bmatrix}\text{prox}_{\iota_{\{|\cdot-\gamma_1|\le \delta\}}} (u_1) &  \text{prox}_{\iota_{\{|\cdot-\gamma_2|\le \delta\}}} (u_2) & \cdots &\text{prox}_{\iota_{\{|\cdot-\gamma_d|\le \delta\}}} (u_{2p}) \end{bmatrix}^\top,
$$
where for $1 \le i \le 2p$
$$
\text{prox}_{\lambda |\cdot|} (u_i) = \max\{|u_i|-\lambda, 0\} \frac{u_i}{|u_i|}
$$
and
$$
\text{prox}_{\iota_{\{|\cdot-\gamma_i|\le \delta\}}} (u_i)=\gamma_i+\max\{|u_i-\gamma_i|, \delta\} \frac{u_i-\gamma_i}{|u_i-\gamma_i|}
$$

Summarizing the above, one has the following proximity operator based algorithm (POA) for approximating a solution to Model~\eqref{model:dantzig overcomplete}.

\begin{algorithm}\caption{(POA)}\label{alg:matrix-final}
	\textbf{Initial Parameters}: The observations $y$ are known, and the sensing matrix $X$ and dictionary $B$ are used to define the matrix $A$, the vector $\gamma$ and the set $\mathcal{F}$.  A parameter $\alpha>0$ is chosen.

	\textbf{Initial Step}: Initial guess $\tau^0=\tau^{-1}=0$, $c^0 = 0$, $\beta^0=0$; $\lambda=0.999\alpha/\|A\|_2^2$.

	\textbf{Main Iterations}: Generate the sequence $\{(c^k,\tau^k): k \in \mathbb{N}\}$ via the iterative scheme~\eqref{eq:itr} until a stopping criterion is met.

	\textbf{Post-processing}: Use the appropriate post-processing scheme to construct the Dantzig estimator $\hat{c}$ from the final output of the main iterative step, and obtain an approximate separation of the signal components from $\hat{c}$ and $B$.
\end{algorithm}

The main iteration process will ideally terminate once the sequence $\{(c^k,\tau^k)\}$ reaches a stationary point.  In practice, we estimate this by stopping once either of the following are met:
\begin{enumerate}
	\item $\displaystyle \frac{\| D^{-1} B^\top X^\top \left(XBc^k -  y\right)\|_{\infty}}{\max\left\{ \|c^k\|_2,1\right\}} \leq \epsilon$, for some $0 < \epsilon < 1$, or
	\item The support of $c^k$ is stationary for $\eta$ iterations, for some predetermined positive integer $\eta$.
\end{enumerate}

As noted in~\cite{C05thedantzig}, the Dantzig selector tends to slightly underestimate solution values within the support of the true solution. 
Let $(c^{\infty},\tau^{\infty})$ denote the final product of the Main Iterations step in POA, and define the set $\Lambda = \text{supp}\{c^{\infty}\}$.  Denote by $c_{\Lambda}$ the vector whose elements are chosen from $c$ with indices in $\Lambda$, and by $B_{\Lambda}$ the submatrix whose columns are chosen from $B$ with indices in $\Lambda$.
The Dantzig selector is estimated on $\Lambda$ by solving the least squares problem $\hat{c}_{\Lambda} = \argmin\{\|X^\top \left(XB_{\Lambda}c -  y\right)\|_2\}$ and setting $\hat{c}_i = 0$ for $i\notin \Lambda$.

The main contribution to the computational complexity of POA is the `Main Iterations' procedure.  Assume the matrices $A$ and $A^\top$ are computed at the beginning of the algorithm, then recalled for use in Equation~\eqref{eq:itr}.  For each iteration of the `Main Iterations' stage, computing $(\tau^k,c^k)$ via Equation~\eqref{eq:itr} requires $\mathcal{O}(q^2)$ multiplications, and determining whether the stopping criteria have been met contributes an additional $\mathcal{O}(q^2)$, where $q$ is the length of the coefficient vector $c$.  Therefore, if $R$ iterations are required to complete POA, then the overall complexity of the algorithm is $\mathcal{O}(Rq^2)$.

\section{Numerical Experiments}
In this section, the separation of noisy undersampled composite signals using POA to solve Model~\eqref{model:dantzig overcomplete} is demonstrated.
All codes are implemented in MATLAB R2013b on a workstation with an Intel i7-3630QM CPU (2.40GHz) and 16GB RAM.

Four numerical experiments will be presented.  In the first three experiments, composite signals are simulated, then POA is applied to the noisy incomplete observations $y = XBc+z$ to recover and separate the components.  The entries of the $n\times p$ sensing matrix $X$ are sampled from the standard normal distribution, which is then normalized so that each column has unit $\ell_2$ norm.  The noise vector $z$ has i.i.d.\ entries sampled from the normal distribution $N(0,\sigma^2)$, with $\sigma = 0.01$ and $0.05$ corresponding to $1\%$ and $5\%$ noise levels respectively.  To measure the effeciency and accuracy of the algorithm, we use the CPU run time and the relative $\ell_2$ error of the recovery of the separate components.  For each experiment, we report the means and standard deviations of these measurements over 50 simulations for each set of parameters.

In the fourth experiment, real-world data taken from the United States Postal Service handwritten digits data sets, are used to sample the composite signal and to train the overcomplete dictionary.  Given an undersampled signal $y = X\beta$, where $X$ is an $n\times p$ sensing matrix with entries taken from a Bernoulli distribution and $\beta$ is the composite image, POA is used to identify the digits and to approximate the individual components in the composite image.

To demonstrate the utility of POA, we also use the Alternating Direction Method (ADM) developed in ~\cite{SPG2, ADM, augmentedLagrangian} to estimate a solution to Model~\eqref{model:dantzig overcomplete} in Experiment 1.  A comparison of the results are summarized in Figure~\ref{fig:experiment 1 time comparison} and Table~\ref{tab:experiment 1 accuracy comparison}.  POA and the ADM separated the components of a composite signal with similar accuracy, however POA was significantly faster.  We do not compare POA and the ADM in Experiments 2 and 3 because these problem domains involve complex-valued dictionaries and coefficients, but the ADM was designed to solve real-valued convex optimization problems only.

The effectiveness of the postprocessing scheme is illustrated in Figure~\ref{fig:experiment 2}.  Every other figure and table present only results with postprocessing.  In the following, time is measured in seconds, and the symbol $\hat{ }$ over a parameter indicates its recovered value via POA with postprocessing.  For notational convenience, the relative $\ell_2$ error of the recovery of a vector $x$ is denoted by $E(x) := \| x - \hat{x}\|_2 / \| x \|_2$.

\subsection{Experiment 1}
In this experiment, we recover and separate a composite signal that has a sparse representation in terms of wavelets and discrete cosine transforms.
Let $\beta$ represent the composite signal, $\beta_\Phi$ the wavelet component and $\beta_\Psi$ the sinusoidal component, each of length $p$.  The signals are simulated by generating a sparse coefficient vector $c$ of length $2p$ by selecting a support set of size $s$ uniformly at random and sampling $c$ on the support set from the distribution $N(100,15)$.  The signals are observed via $y = XBc + z$, where the overcomplete dictionary $B$ is formed by concatenating matrices for the level 5 Haar wavelet decomposition and the discrete cosine transform.  The simulations are performed 50 times for each $m = 1,2,\ldots,10$ and each $\sigma = 0.01, 0.05$, with parameters $p = 256m+512,\; n = p/4,\; s = \lceil n/9 \rceil$.  The parameters for the stopping criteria of POA are $\epsilon = 10^{-4}$ and $\eta=20$.

\begin{figure}[htbp]
	\centering
	\begin{subfigure}[t]{.48\textwidth}
		\includegraphics[width=\textwidth, height=.25\textheight]{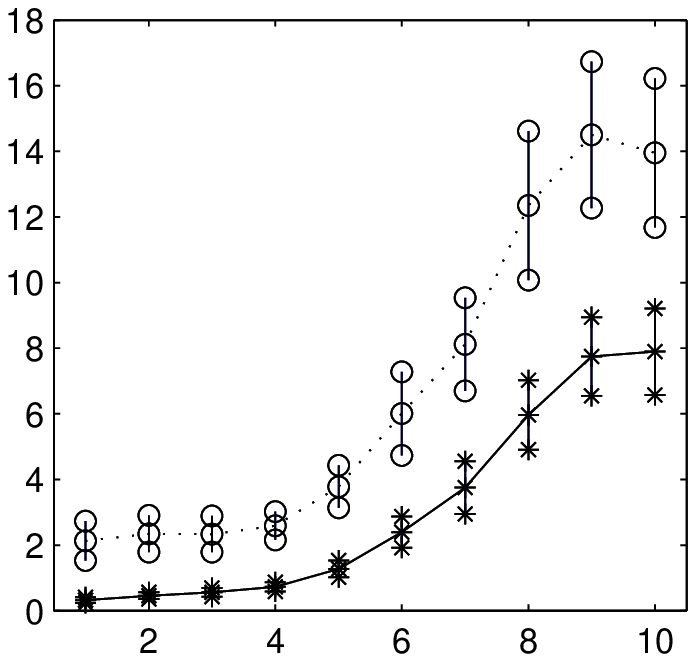}
		\caption{CPU runtime, with $\sigma = 0.01$.}
		\label{subfig:CPU low}
	\end{subfigure}
	\begin{subfigure}[t]{.48\textwidth}
		\includegraphics[width=\textwidth, height=.25\textheight]{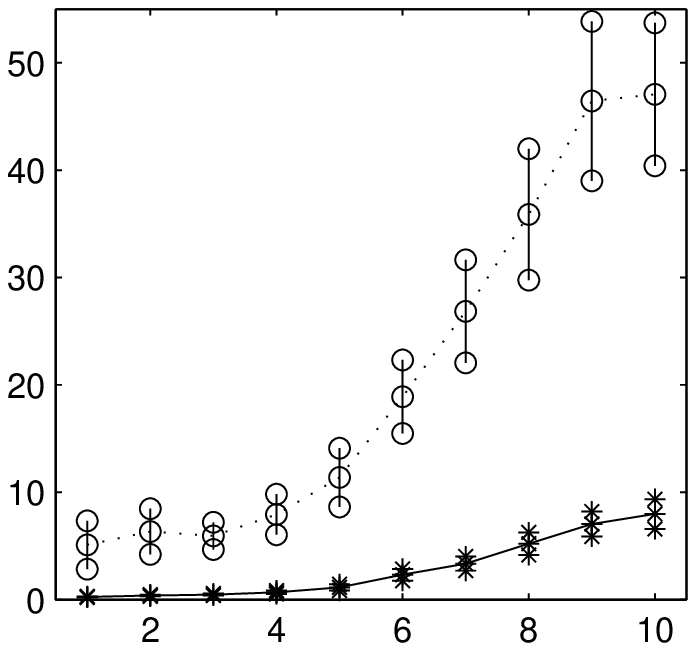}
		\caption{CPU runtime, with $\sigma = 0.05$.}
		\label{subfig:CPU high}
	\end{subfigure}
	\caption{Comparisons of CPU runtime of POA and ADM from Experiment 1.  The mean CPU runtime of the 50 simulations of each algorithm is represented by the points along the curves, and the vertical lines represent one standard deviation from the means.  The solid line with * markers plot the values obtained using POA, and the dotted line with $\circ$ markers plot the values obtained using the ADM.  The horizontal axis represents the parameter $m$ determining the size of the system, and the vertical axis is measured in seconds.}
	\label{fig:experiment 1 time comparison}
\end{figure}

Comparisons of the mean and standard deviation of CPU time and relative $\ell_2$ recovery error of each of the components for 50 simulations at each level of $m$ and $\sigma$  using both POA and the ADM to solve Model~\eqref{model:dantzig overcomplete} are illustrated in Figure~\ref{fig:experiment 1 time comparison} and Table~\ref{tab:experiment 1 accuracy comparison}. \hfill $\Box$

\begin{table}
	\centering
	\begin{tabular}{c c||c c|c c||c c|c c}
	&$m$
	&\multicolumn{4}{c||}{$E(\beta_\Phi)$}
	&\multicolumn{4}{c}{$E(\beta_\Psi)$}\\
	&& \multicolumn{2}{c}{POA}
	&\multicolumn{2}{c||}{ADM}
	&\multicolumn{2}{c}{POA}
	&\multicolumn{2}{c}{ADM}\\
	\hline \hline
	&&Mean	&Std Dev	&Mean	&Std Dev	&Mean 	&Std Dev	&Mean	&Std Dev\\
	&&$(\times 10^{-3})$ &$(\times 10^{-3})$ &$(\times 10^{-3})$ &$(\times 10^{-3})$ &$(\times 10^{-3})$ &$(\times 10^{-3})$ &$(\times 10^{-3})$ &$(\times 10^{-3})$\\	
	\hline
	\multirow{10}{*}{\begin{sideways}Noise Level $\sigma = 0.01$\end{sideways}}
	&1	&9.0421	&2.2306		&8.6562	&2.0236	&9.1171	&2.9345	&8.7136	&2.7528	\\
   	&2	&8.3235	&2.2061		&8.1769	&2.2004	&8.3946	&2.0634	&8.3330	&2.2880	\\
   	&3	&8.3875	&2.2645		&8.3397	&2.2484	&8.6408	&2.4491	&8.5115	&2.3044	\\
   	&4	&7.7607	&1.3936		&7.7985	&1.3568	&7.6132	&1.5551	&7.6046	&1.5981	\\
   	&5	&8.2638	&1.4670		&8.2121	&1.6065	&7.5479	&1.3118	&7.6928	&1.4544	\\
   	&6	&8.0279	&1.4589		&8.0491	&1.5101	&7.7504	&1.4392	&7.7910	&1.6199	\\
   	&7	&7.7411	&1.2514		&7.6378	&1.1422	&7.9663	&1.4901	&7.9039	&1.4730	\\
   	&8	&8.2142	&1.3920		&8.1150	&1.3064	&7.9014	&1.3367	&7.8035	&1.3457	\\
   	&9	&7.8209	&1.0384		&7.8440	&1.1215	&8.1155	&1.4212	&8.0859	&1.4590	\\
   	&10	&7.9987	&1.2901		&7.9806	&1.2391	&7.9123	&1.2342	&7.9213	&1.1699	\\
	\hline		
	&&Mean	&Std Dev	&Mean	&Std Dev	&Mean 	&Std Dev	&Mean	&Std Dev\\
	&&$(\times 10^{-2})$	&$(\times 10^{-3})$
	&$(\times 10^{-2})$ 		&$(\times 10^{-3})$
	&$(\times 10^{-2})$ 		&$(\times 10^{-3})$
	&$(\times 10^{-2})$ 		&$(\times 10^{-3})$\\	
	\hline
	\multirow{10}{*}{\begin{sideways}Noise Level $\sigma=0.05$\end{sideways}}
	&1	&4.4187	&14.921		&4.2647	&13.274	&4.2368	&13.321	&4.1104	&12.031	\\
	&2	&3.9471	&10.951		&3.8098	&10.353	&3.9301	&10.570	&3.9203	&10.816	\\
	&3	&3.9268	&9.6222		&3.8574	&9.0583	&3.9068	&7.8635	&3.8550	&7.4182	\\
	&4	&3.9723	&7.8517		&3.9302	&7.3932	&4.0431	&12.237	&3.9847	&11.224	\\
	&5	&4.0825	&8.4696		&4.0847	&8.1707	&3.8121	&6.9432	&3.8310	&6.7285	\\
	&6	&3.9041	&4.9101		&3.8963	&5.0425	&3.8486	&7.1261	&3.7879	&7.3865	\\
	&7	&4.0000	&5.8845		&3.9813	&5.4931	&4.1028	&8.0082	&4.1050	&7.6549	\\
	&8	&3.8323	&6.2478		&3.8310	&6.3490	&3.9243	&7.4910	&3.9327	&7.7148	\\
	&9	&3.8040	&6.4910		&3.8195	&6.3716	&3.9178	&5.7576	&3.9010	&5.6531	\\
	&10	&3.9514	&6.2999		&3.9518	&6.2383	&3.9898	&6.3100	&3.9677	&6.3825	\\
	\hline
	\end{tabular}
	\caption{Comparisons of the relative $\ell_2$ errors of the recovered components from Experiment 1 using POA and ADM. The means and standard deviations over 50 simulations of the relative errors are given for the recovery of each component in the original signal, and for each parameter $m$ and noise level $\sigma$.}
	\label{tab:experiment 1 accuracy comparison}
\end{table}

\subsection{Experiment 2}
In this experiment, we recover and separate a composite signal, with one component $\beta_\Phi$ being a dirac spike signal with random sparse locations and values, and the other component $\beta_{\Psi}$ being a sinusoidal signal with random sparse Fourier coefficients.  The parameters used in this experiment are $p = 256m,\; n=64m,$ and $s=\lceil n/9\rceil$ for $m=1,2,\ldots,10$, and stopping criteria parameters $\epsilon=10^{-4}$ and $\eta=6$ for $\sigma=0.01$ and $\eta=30$ for $\sigma=0.05$.

For each of the 50 simulations for each value of $m$ and $\sigma$, the experiment is set up by generating a vector $c$ of length $2p$ by selecting a support set $S$ of size $s$ uniformly at random, then sampling $c$ on $S$ with i.i.d. entries $C(j) = \lambda(j)(1+|a(j)|)$, where $\lambda(j)$ is $1$ or $-1$ and $a(j)\sim N(0,1)$.  The signal is observed via the vector $y = XBc + z$, with the dictionary $B$ being a concatenation of the identity matrix and the discrete Fourier transform matrix, both of size $p\times p$.

\begin{figure}
	\centering
	\begin{subfigure}[t]{0.48\textwidth}
		\includegraphics[width=\textwidth]{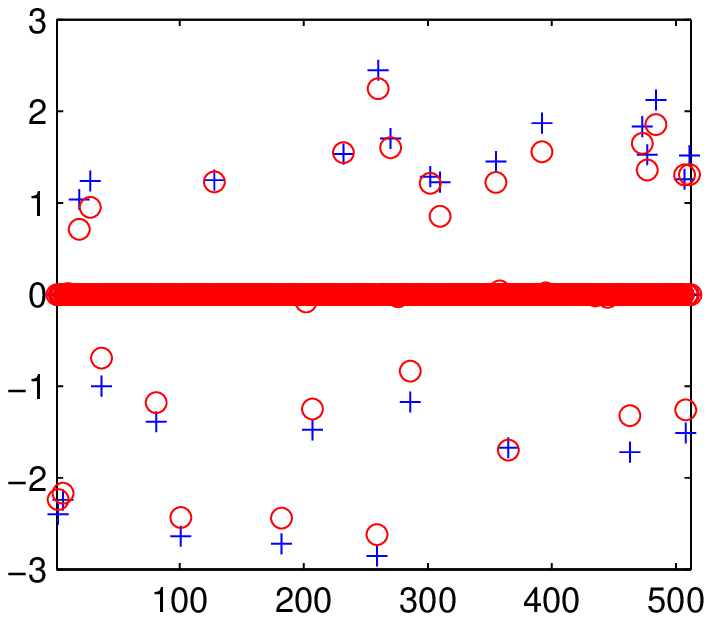}
		\caption{Recovery of $\beta_\Phi$ without postprocessing.}
	\end{subfigure}
	\begin{subfigure}[t]{0.48\textwidth}
		\includegraphics[width=\textwidth]{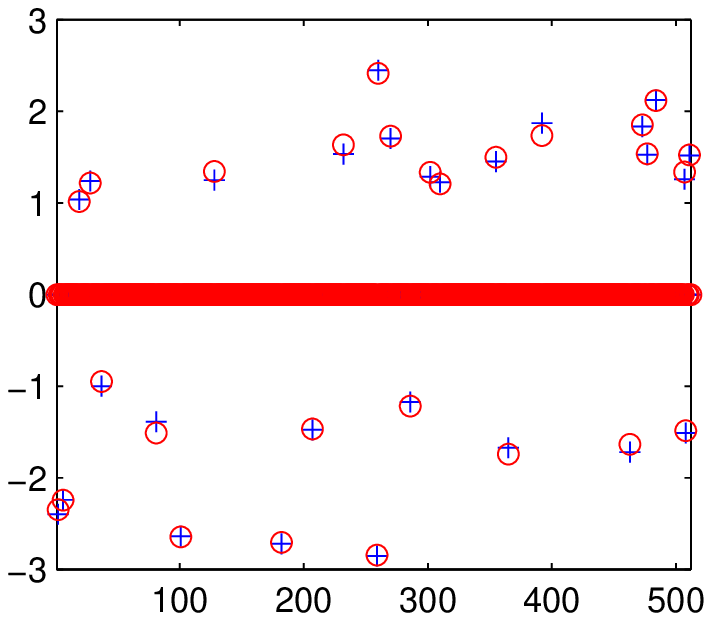}
		\caption{Recovery of $\beta_\Phi$ with postprocessing.}
	\end{subfigure}

	\begin{subfigure}[t]{0.48\textwidth}
		\includegraphics[width=\textwidth]{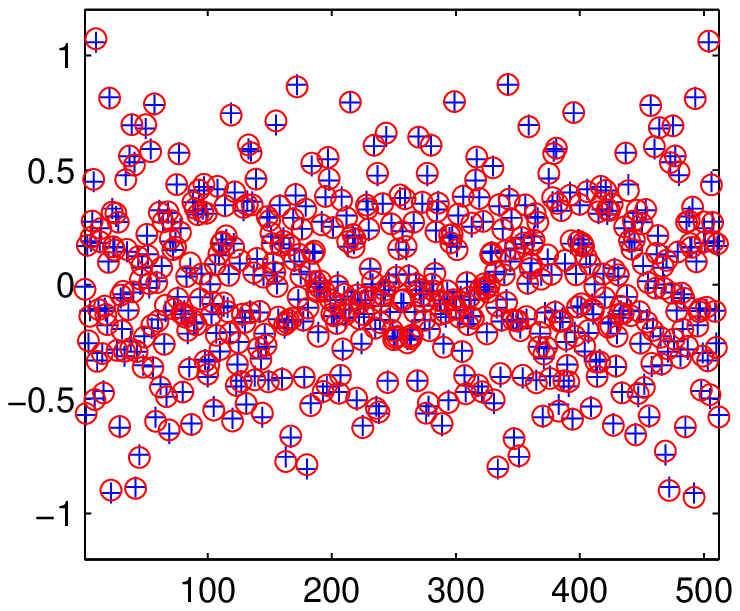}
		\caption{Recovery of $\mathfrak{Re}(\beta_\Psi)$.}
	\end{subfigure}
	\begin{subfigure}[t]{0.47\textwidth}
		\includegraphics[width=\textwidth]{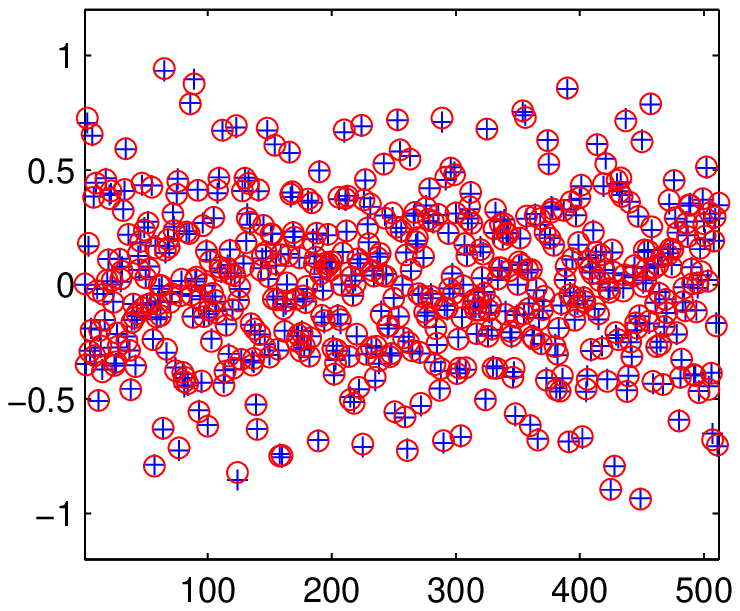}
		\caption{Recovery of $\mathfrak{Im}(\beta_\Psi)$.}
	\end{subfigure}
\caption{The separation of components from a typical simulation of Experiment 2 with $m = 2, \sigma = 0.05$.  The plots in the top row are of the recovered values of the component $\beta_\Phi$ without and with postprocessing.  The component $\beta_\Psi$ is complex-valued, so the recovery of its real and imaginary parts with postprocessing are displayed separately in the second row.  In each plot, the exact values are denoted by `+' and the values recovered by POA are denoted by `o'.
}\label{fig:experiment 2}
\end{figure}

Figure~\ref{fig:experiment 2} illustrates the accuracy of the numerical separation by comparing the numerically recovered composite signal and separated components (denoted by `o') against the exact values of the signal and components (denoted by `+') for one simulation with $m=2$ and $\sigma=0.05$. Table~\ref{tab:experiment 2} lists the mean and standard deviation of the algorithm run time and the relative $\ell_2$ recovery error of each signal component over the 50 simulations for each $m$ and $\sigma$. 
\hfill $\Box$

\begin{table}
	\centering
	\begin{tabular}{c c||c c|c c|c c}
	&$m$		
	&\multicolumn{2}{c|}{Time (Seconds)} 	
	&\multicolumn{2}{c|}{$E(\beta_\Phi)$}
	&\multicolumn{2}{c }{$E(\beta_\Psi)$}\\
	\hline
	\hline
	&&\multirow{2}{*}{Mean}	&\multirow{2}{*}{Std Dev}	 &Mean	&Std Dev &Mean	&Std Dev\\
	&&&&$(\times 10^{-3})$	&$(\times 10^{-4})$		&$(\times 10^{-3})$	&$(\times 10^{-4})$\\
	\hline
	\multirow{10}{*}{\begin{sideways}Noise Level $\sigma = 0.01$\end{sideways}}
   &1		&2.8118$\times 10^{-2}$  	&5.8329$\times 10^{-3}$	&6.4762	&15.851    &6.7300	&12.374	\\
   &2 	&6.2733$\times 10^{-2}$  	&2.1378$\times 10^{-2}$	&6.4622	&8.6297	&6.7490	&8.6886	\\
   &3 	&0.1351		   		&7.5995$\times 10^{-2}$	&6.6026	&8.7737	&6.4949	&7.8729	\\
   &4 	&0.4025		   		&0.2948				&6.7021	&8.2956	&6.6030	&6.3652	\\
   &5		&0.9603   				&0.7242 				&6.7423	&6.6547	&6.7797	&6.6109	\\
   &6		&1.5431	     			&1.3161				&6.6408	&5.1202	&6.6585	&4.4381	\\
   &7		&3.0107         			&2.2266				&6.8240	&5.9514	&6.8955	&6.7195	\\
   &8 	&3.5928         			&2.6642				&6.8223	&6.2282	&6.6457	&4.2048	\\
   &9		&5.2908				&3.7904				&6.6667	&5.3515	&6.5987	&4.0810	\\
   &10	&8.6188				&4.7235				&6.8529	&5.7726	&6.8113	&4.5567	\\
	\hline
	&&\multirow{2}{*}{Mean}	&\multirow{2}{*}{Std Dev}	 &Mean	&Std Dev &Mean	&Std Dev\\
	&&&&$(\times 10^{-2})$	&$(\times 10^{-3})$		&$(\times 10^{-2})$	&$(\times 10^{-3})$\\
	\hline
	\multirow{10}{*}{\begin{sideways}Noise Level $\sigma = 0.05$\end{sideways}}
   &1		&0.1076	   	&3.3335$\times 10^{-2}$  	&3.4973	&6.6099	&3.3185	&6.6436	\\
   &2		&0.2615	   	&6.4574$\times 10^{-2}$	&3.5430	&4.2359	&3.4794	&4.3713	\\
   &3		&0.6093	   	&0.1218	   			&3.4250	&4.4838	&3.3384	&3.4251	\\
   &4		&1.7032	      	&0.3350	   			&3.3804	&3.2731	&3.3572	&3.2625	\\
   &5		&3.6084          	&0.5986	   			&3.4346	&2.8969	&3.3914	&2.6619	\\
   &6		&5.5694 		&0.9674	   			&3.4677	&2.8924	&3.3116	&2.6373	\\
   &7		&7.5435 		&1.1237				&3.4493	&2.7535	&3.3879	&1.9975	\\
   &8		&9.9614 		&1.1852				&3.4638	&2.3516	&3.4088	&2.2992	\\
   &9		&12.8380		&1.7008				&3.4874	&2.3352	&3.3968	&2.5498	\\
   &10	&15.4924 		&1.4241				&3.4326	&2.3453	&3.4374	&1.4001	\\
	\hline
	\end{tabular}
	\caption{The means and standard deviations of the CPU run time and the relative $\ell_2$ recovery error of each component for 50 simulations of Experiment 2 using POA for each parameter $m$ and noise level $\sigma$.}
	\label{tab:experiment 2}
\end{table}

\subsection{Experiment 3}
In this experiment, a specified composite signal is separated into distinct components.  The composite signal $\beta$ is formed by combining the discretized sinusoidal component $\beta_\Phi$ and the sparse spike component $\beta_\Psi$, where
\[ \beta_{\Phi}(x) = 30\sin\left(\frac{2\pi x}{p}\right) + \sin\left(\frac{\pi x}{2}\right), \text{ for } x \in \{0,1,\ldots,1023\},\]
and $\beta_{\Psi}$ is formed by selecting a set $S$ of size $s$ and sampling on $S$ using the same method applied to $c$ in Experiment 2.
The overcomplete dictionary $B$ is the concatenation of the identity and the discrete Fourier transform matrices, each of size $p\times p$.  Note that in this experiment, the exact coefficient vector is not directly specified, so the signal is observed according to $y = X\beta + z$.  The parameters used are $p = 1024,\; n = 512,\; s = 57$, with stopping criteria $\epsilon = 10^{-6}$ and $\eta=6$ for $\sigma =0.01$ and $\eta=30$ for $\sigma = 0.05$.

\begin{figure}[htbp]
	\centering
	\begin{subfigure}[t]{0.40\textwidth}
		\includegraphics[width=\textwidth]{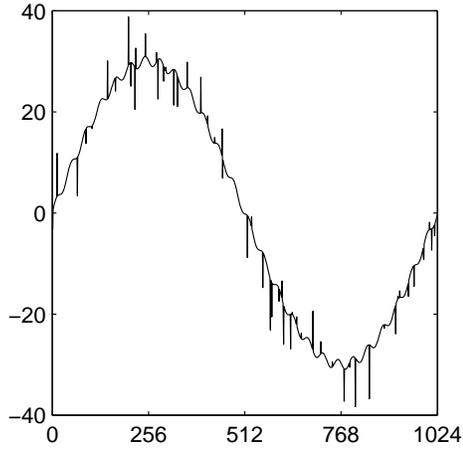}
		\caption{The exact signal $\beta$.}
	\end{subfigure}
	\begin{subfigure}[t]{0.40\textwidth}
		\includegraphics[width=\textwidth]{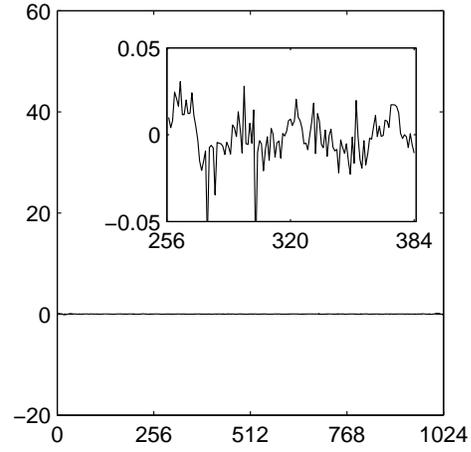}
		\caption{Recovery error $\beta - \hat{\beta}$.}
	\end{subfigure}

	\begin{subfigure}[t]{0.40\textwidth}
		\includegraphics[width=\textwidth]{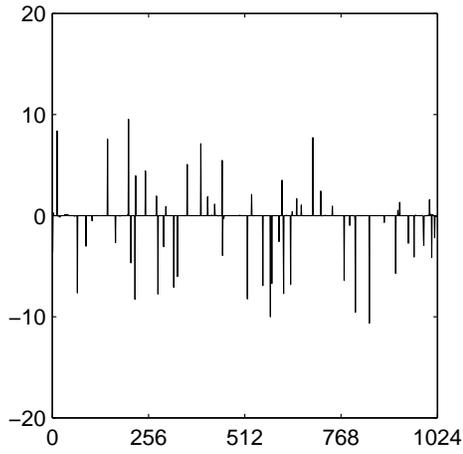}
		\caption{The recovered component $\mathfrak{Re}(\hat{\beta}_\Phi)$.}
	\end{subfigure}
	\begin{subfigure}[t]{0.40\textwidth}
		\includegraphics[width=\textwidth]{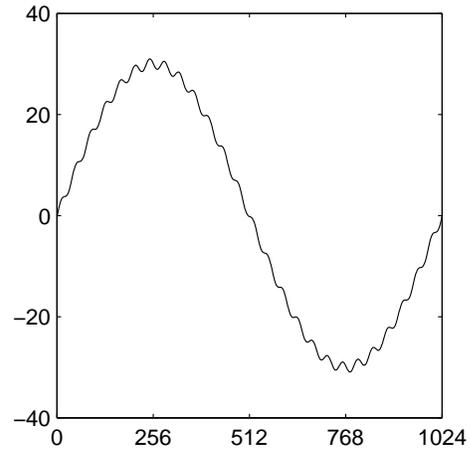}
		\caption{The recovered component $\hat{\beta}_\Psi$.}
	\end{subfigure}
	\caption{The separation of components from a typical simulation of Experiment 3 with noise level $\sigma = 0.01$.  Subplot (a) is the exact composite signal $\beta$, which is then undersampled via $y = X\beta + z$.  Subplot (b) is the pointwise error $\beta - \hat{\beta}$, where $\hat{\beta}$ is the recovered signal.  Subplots (c) and (d) are the individual recovered components $\mathfrak{Re}(\hat{\beta}_\Phi)$ and $\hat{\beta}_\Psi$ separated using POA and Model~\eqref{model:dantzig overcomplete}.}
	\label{fig:experiment 3}
\end{figure}

A typical plot of the composite signal, the individual components, their recoveries and pointwise recovery error with $\sigma = 0.01$ are shown in Figure~\ref{fig:experiment 3}. Table~\ref{tab:experiment 3} displays the means and standard deviations of the run time of POA and relative $\ell_2$ recovery error of the signal and each component over 50 simulations for each noise level. \hfill $\Box$

\begin{table}
	\centering
	\begin{tabular}{c|c c|c c}
	&\multicolumn{2}{c|}{$\sigma = 0.01$}	&\multicolumn{2}{c}{$\sigma = 0.05$}\\
	&Mean	&Std Dev	&Mean	&Std Dev\\
	\hline
	Time		
	&2.9154	&5.7204$\times 10^{-2}$ 	&2.8984	&8.6182$\times 10^{-2}$\\
	$E(\beta)$ 
	&2.3460$\times 10^{-3}$	&2.0277$\times 10^{-3}$ 	&5.3163$\times 10^{-3}$	&3.2198$\times 10^{-3}$\\
	$E(\beta_\Phi)$ 
	&2.0391$\times 10^{-3}$ 	&1.5152$\times 10^{-3}$	&3.7928$\times 10^{-3}$	&2.3483$\times 10^{-3}$\\
	$E(\beta_\Psi)$ 
	&3.8823$\times 10^{-2}$	&3.3788$\times 10^{-2}$	&8.4469$\times 10^{-2}$	&5.9878$\times 10^{-2}$\\
	\hline
	\end{tabular}
	\caption{The means and standard deviations of the run time of POA and the relative $\ell_2$ recovery norm of the composite signal and each component over 50 simulations of Experiment 3 for noise levels $\sigma = 0.01$ and $0.05$.  }
	\label{tab:experiment 3}
\end{table}

\subsection{Experiment 4}
In this experiment, composite signals of two handwritten digits are classified and separated using POA and standard principal component analysis techniques.  The data are taken from the USPS handwritten digit data sets, obtained from~\cite{digits}.  The data set contains 10 classes, labeled zero through nine, of 1100 examples of 8-bit grayscale $16\times 16$ images.  Each image is reshaped as $256$-column vectors.  In each class, 998 examples are used as the training set and the remaining 102 examples form the testing set.  Denote by $R[j]$ the collection of vectors forming the training set, and by $T[j]$ the collection of vectors forming the testing set.  
By an abuse of notation, also denote by $R[j]$ the matrix whose columns come from the $j^\text{th}$ collection of training vectors.  

The composite vector $\beta$ is sampled from the test data by randomly selecting two vectors from two distinct test sets, $\beta_1 \in T_{j_1}$ and $\beta_2 \in T_{j_2}$, and taking $\beta = \beta_1 + \beta_2$.  Consider the observation $y = X\beta$, where $X$ is a random Bernoulli $128 \times 256$ sensing matrix.  The overcomplete dictionary is learned from the labeled training sets by concatenating the first few principal components of each matrix $R[j]$.  Using Model~\eqref{model:dantzig overcomplete} and POA, we classify the two digits in the composition $\beta$, and recover an approximation of the two digits using the procedure outlined in Algorithm~\ref{alg:handwritten-digits}.

\begin{algorithm}\caption{Composite Handwritten Digit Classification and Separation}\label{alg:handwritten-digits}
	\textbf{(1)}: Input the observation $y$, the training sets $R[0], \ldots, R[9]$, and a positive integer $k\leq256$.

	\textbf{(2)}: For each $j$, compute $\tilde{U}[j]$, a matrix whose columns are the first $k$ principal components of the matrix $R[j]$.  The overcomplete dictionary is learned from the training data via 
	\[ B =  \begin{bmatrix} \tilde{U}[0] & \tilde{U}[1] & \cdots & \tilde{U}[9] \end{bmatrix}.\]

	\textbf{(3)}: Apply POA to $y$ to produce a sparse coefficient vector $\hat{c}$ satisfying Model~\eqref{model:dantzig overcomplete} and to recover the composite vector $\hat{\beta} = B\hat{c}$.

	\textbf{(4)}: Classify the components of $\hat{\beta}$ as the indices $j_1, j_2$ yielding the smallest two values of 
	\[ \left\| \left( I_{256} - \tilde{U}[j]\tilde{U}[j]^* \right) \hat{\beta}\right\|_2.\]

	\textbf{(5)}:  Form a reduced dictionary 
	\[ \hat{B} = \begin{bmatrix} \tilde{U}[j_1] & \tilde{U}[j_2] \end{bmatrix}.\]
	
	\textbf{(6)}: Apply POA again to $y$ using $\hat{B}$ to obtain a new coefficient vector $\hat{c} = \begin{bmatrix} c_{j_1}^\top & c_{j_2}^\top \end{bmatrix}^\top$ satisfying Model~\eqref{model:dantzig overcomplete}.  The two unknown components of $\beta$ are approximated as 
	\[ \hat{\beta}_{j_1} = \tilde{U}[j_1]\hat{c}_{j_1},\; \text{ and } \; \hat{\beta}_{j_2} = \tilde{U}[j_2]\hat{c}_{j_2}.\]
\end{algorithm}

Let us explain the motivation behind Algorithm~\ref{alg:handwritten-digits}.  Since each unknown digit in the composite vector should be described well by its corresponding principal vectors, it follows that one can approximate the composite vector by $\beta = Bc$ for a sparse coefficient vector $c$, and Model~\eqref{model:dantzig overcomplete} is appropriate to use on $y$.  In Step (2), the first $k$ principal components of each training set $R[j]$ is computed and used to form the overcomplete dictionary.  One possible method to find each $\tilde{U}[j]$ is to use the singular value decomposition.  If $ U\Sigma V^*$ is a singular value decomposition of the $R[j]$,  then $\tilde{U}[j]$ can be taken as the first $k$ columns of the collection of left singular vectors $U$. In Step (4), the recovered composite vector $\hat{\beta}$ is projected onto the vector spaces spanned by the first $k$ principal components of each training class, and the components are classified according to the residual vectors giving the smallest $\ell_2$ norm.  In Steps (5) and (6), the dictionary is reduced based on the most significant principal components identified in Step (4).  Reducing the dictionary results in a cleaner separation since only coefficients in the identified classes will be recovered.  

\begin{figure}[thbp]
	\centering
	\begin{subfigure}[t]{0.48\textwidth}
		\includegraphics[width=\textwidth]{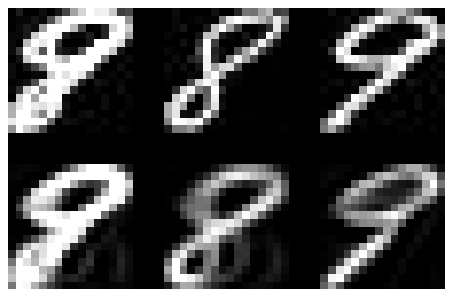}
	\end{subfigure}
	\begin{subfigure}[t]{0.48\textwidth}
		\includegraphics[width=\textwidth]{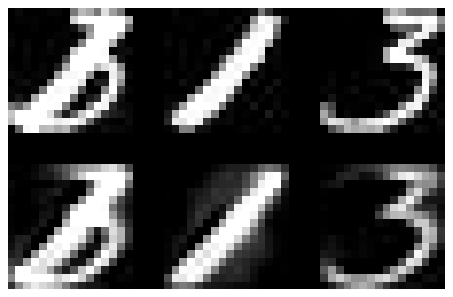}
	\end{subfigure}
	\caption{Two illustrations of the approximation of individual components composite images using POA in Algorithm~\ref{alg:handwritten-digits} with an overcomplete dictionary learned from training examples as in Experiment~4.  In each collection, reading left to right, the top row is the exact composite image $\beta$ followed by the exact components.  The bottom row is the recovered composite image $\hat{\beta}$ followed by the recovered components $\hat{\beta}_{j_1}$ and $\hat{\beta}_{j_2}$.}
	\label{fig:experiment 4}
\end{figure}

Algorithm~\ref{alg:handwritten-digits} was performed 1000 times for different observations $y$ with $k=30$, and the recovered classification indices were compared with the true classification of the two components, as well as the classification determined using the smallest residual error from the composite image $\beta$.  In 95.4\% of the simulations, the classification accuracy of $\hat{\beta}$ generated by POA given only the random projections $y$ matched or exceeded the classification accuracy of the exact image $\beta$.  Examples of the separation of two composite images sampled from the testing data sets using Algorithm~\ref{alg:handwritten-digits} is illustrated in Figure~\ref{fig:experiment 4}.  \hfill $\Box$

\section{Conclusion}

The strength of Model~\eqref{model:dantzig overcomplete} and POA in separating noisy undersampled composite signals is demonstrated through the numerical experiments in the previous section.  Figures~\ref{fig:experiment 2}, \ref{fig:experiment 3} and~\ref{fig:experiment 4} and Tables~\ref{tab:experiment 1 accuracy comparison}, \ref{tab:experiment 2} and~\ref{tab:experiment 3} clearly illustrate components of a composite signal are separated using POA with a high degree of precision.  In particular, Figure~\ref{fig:experiment 2} demonstrates the method's ability to distinguish mixed signals that have similar dynamic range, Figure~\ref{fig:experiment 3} demonstrates the effectiveness of this method to separate a smooth signal from unwanted impulse noise and Figure~\ref{fig:experiment 4} demonstrates POA can be used to separate overlaid images using an overcomplete dictionary training from labeled data.

The relative $\ell_2$ norm error of the recovery of each of the components remains fairly stables as $m$ increases for each $\sigma$, as shown in Tables~\ref{tab:experiment 1 accuracy comparison} and~\ref{tab:experiment 2}. This suggests that the method is not greatly affected by the system size, and is stable and reliable in practice.

Comparisons to ADM demonstrate the advantages of using POA to separate a noisy undersampled composite signal and estimate solutions of Model~\eqref{model:dantzig overcomplete}.  When the underlying coefficient vector and overcomplete dictionary have only real-valued elements, the two methods have similar accuracy yet POA is significantly faster (Figure~\ref{fig:experiment 1 time comparison}, Table~\ref{tab:experiment 1 accuracy comparison}).  Moreover, POA can be readily applied to models involving complex-valued overcomplete dictionaries and coefficients, but the ADM used for comparison was designed for convex optmization in real-valued domains only.  In this sense, POA applies to a wider class of problems.

In summary, we have introduced a model for the Dantzig selector incorporating overcomplete dictionaries that can be used to separate composite signals. Additionally, we have proposed POA, an iterative algorithm to estimate solutions to Model~\eqref{model:dantzig overcomplete} and have given the results of several numerical experiments to support the strength of both the model and the algorithm.  We have shown through the numerical experiments that POA is preferable over the competing method ADM, since POA produces results with similar accuracy but with less CPU time and is applicable to a wider class of problems.   Some possible future applications based on the foundation of the separation of composite signals include medical imaging, image inpainting and feature extraction.  Moreover, advances in finding sparse descriptions of noisy composite signals can be applied to improve technologies in signal compression and communications.

\bibliographystyle{acm}
\bibliography{CSbib}

\end{document}